# Вариации на тему уточнённого порядка

Б. Н. Хабибуллин

*Башкирский государственный университет*
*Россия, Республика Башкортостан, 450076 г. Уфа, улица Заки Валиди, 32.*

*Email: khabib-bulat@mail.ru*

Понятие уточнённого порядка широко используется в теориях целых, мероморфных, субгармонических и плюрисубгармонических функций. Приводится общая трактовка этого понятия как уточнённой функции роста относительно модельной функции роста. Классический уточнённый порядок – это $\ln V$, когда V – уточнённая функция роста относительно тождественной функции на положительной полуоси. Наше определение использует лишь одно условии. Такая форма определения новая и для классического уточнённого порядка.

Понятие уточнённого порядка возникло в работах Ж. Валирона [1; Ch. III, I.6]. Основные его свойства и применения изложены в [2; гл. I, § 12], [3; гл. II, § 2], [4; 7.4], [5]. Перейдём к точным определениям и формулировкам. $\mathbb{R}$ – множество *действительных* чисел; $\mathbb{R}^+ := \{r \in \mathbb{R}: r \geq 0\}$ – множество *положительных* чисел; $\mathbb{R}^+_* := \mathbb{R}^+ \setminus \{0\}$ –множество *строго положительных* чисел. Подмножество $R \subset \mathbb{R}$ – *луч положительного направления на* $\mathbb{R}$, если $\emptyset \neq R \neq \mathbb{R}$, а также для любых $r \in R$ и $x \in \mathbb{R}$ из $r \leq x$ следует, что $x \in R$. Далее через $R_{\rightarrow}$ обозначаем *какой-либо произвольный луч положительного направления на* $\mathbb{R}$. Рассматриваются только функции $f$ с *областью значений* $\mathbb{R}$ и с *областью определения*, включающей в себя $R_{\rightarrow}$ . Функция $f$ обладает некоторым свойством, если она обладает им на $R_{\rightarrow}$ . Функция $f$ *положительная* и пишем $f \geq 0$ (соответственно *строго положительная* и пишем $f > 0$), если $f(R_{\rightarrow}) \subset \mathbb{R}^+$ (соответственно $f(R_{\rightarrow}) \subset \mathbb{R}^+_*$). Функция $f$ *возрастающая* (соответственно *строго возрастающая*), если для любых $x_1, x_2 \in R_{\rightarrow}$ из $x_1 < x_2$ следует, что $f(x_1) \leq f(x_2)$ (соответственно $f(x_1) < f(x_2)$). Действия над функциями производятся поточечно на $R_{\rightarrow}$ .

Функция $M$ *выпукла относительно* $\ln$, если функция $m: x \mapsto M(e^x)$ *выпуклая.* Из известных свойств выпуклых функций это означает, что функция $M$ *непрерывна* и существуют *левые* и *правые производные* соответственно $M'_-$ и $M'_+$, для которых функция $r \mapsto rM'_-(r)$ и/или функция $r \mapsto rM'_+(r)$ возрастающая [6, Ch. I]. Другими словами, по переменной $z \in \mathbb{C}$ на *комплексной плоскости* $\mathbb{C}$ функция $z \mapsto M(|z|)$ *субгармоническая радиальная функция* вне некоторого круга с центром в нуле из $\mathbb{C}$ [6, Ch. III].

**Определение.** Выпуклую относительно $\ln$ функцию $M > 0$ с производной $M' > 0$ и с пределом $\lim\limits_{r\to+\infty} M(r) = +\infty$ называем *модельной функцией роста*. Дифференцируемую функцию $V > 0$ называем *уточнённой функцией роста относительно модельной функции роста* $M$, если существует хотя бы один из пределов

$$\lim_{r\to+\infty} \frac{M(r)V'(r)}{M'(r)V(r)} = \lim_{r\to+\infty} \frac{(\ln V(r))'}{(\ln M(r))'} = \lim_{x\to+\infty} \frac{(\ln v(x))'}{(\ln m(x))'} = \lim_{x\to+\infty} \frac{m(x)v'(x)}{m'(x)v(x)} \in \mathbb{R}^+, \qquad (1)$$

где, как и выше, функция $m: x \mapsto M(e^x)$ выпукла, а $v: x \mapsto V(e^x)$ дифференцируема, а равенства в (1) при условии существования хотя бы одного из пределов в (1), как нетрудно показать, всегда имеют место.

Использование в (1) вместе с $M$ и $V$ функций $m(\ln r) = M(r)$ и $v(\ln r) = V(r)$ для $r \in \mathbb{R}_*^+$ обусловлено необходимостью вписать в будущем наш подход к функциям роста в общую идеологию по этому направлению, разработанную К. Кизельманом в [7].

**Теорема.** *Пусть $M$ – модельная функция роста, $V > 0$ – дифференцируемая функция. Тогда эквивалентны следующие два утверждения:*

I. *Функция $V$ – уточнённая функция роста относительно модельной функции $M$.*

II. *Для функции*

$$\boldsymbol{\rho}_M := \frac{\ln V}{\ln M} \tag{2}$$

*существуют два конечных предела*

$$\rho := \lim_{r\to+\infty} \boldsymbol{\rho}_M(r) = \lim_{r\to+\infty} \frac{\ln V(r)}{\ln M(r)} \in \mathbb{R}^+, \tag{3}$$

$$\lim_{r\to+\infty} \frac{M(r)}{M'(r)} \boldsymbol{\rho}'_M(r) \ln M(r) = \lim_{r\to+\infty} \frac{\ln M(r)}{(\ln M)'(r)} \boldsymbol{\rho}'_M(r) = 0. \tag{4}$$

*При выполнении любого из этих двух утверждений имеем равенства*

$$\rho := \lim_{r\to+\infty} \boldsymbol{\rho}_M(r) = \lim_{r\to+\infty} \frac{M(r)V'(r)}{M'(r)V(r)}. \tag{5}$$

*Доказательство.* Вычисление производной функции (2) даёт

$$\boldsymbol{\rho}'_M = \frac{(V'/V)\ln M - (M'/M)\ln V}{(\ln M)^2} = \frac{M'}{M\ln M}\left(\frac{MV'}{M'V} - \frac{\ln V}{\ln M}\right),$$

откуда по определению (2) функции $\boldsymbol{\rho}_M$ имеем

$$\frac{MV'}{M'V} = \frac{\ln V}{\ln M} + \frac{M\ln M}{M'} \boldsymbol{\rho}'_M = \boldsymbol{\rho}_M + \frac{M}{M'} \boldsymbol{\rho}'_M \ln M. \tag{6}$$

Отсюда, если выполнено утверждение II Теоремы и существуют пределы (3) и (4), то

$$\lim_{r\to+\infty} \frac{M(r)V'(r)}{M'(r)V(r)} = \lim_{r\to+\infty} \boldsymbol{\rho}_M(r) + \lim_{r\to+\infty} \frac{M(r)}{M'(r)} \boldsymbol{\rho}'_M(r) \ln M(r) = \rho + 0 = \rho \in \mathbb{R}^+.$$

Таким образом, существуют пределы (1), справедливо равенство (5) и по Определению выполнено утверждение I Теоремы.

Обратно, пусть выполнено утверждение I Теоремы, т.е. существует предел

$$\rho := \lim_{r\to+\infty} \frac{(\ln V(r))'}{(\ln M(r))'} \in \mathbb{R}^+. \tag{7}$$

**Правило Лопиталя** [8; теорема 1]**.** *Пусть функции $f$ и $g$ дифференцируемы на некотором луче $R_{\rightarrow}$ положительного направления на $\mathbb{R}$. Если существуют пределы*

$$\lim_{r\to+\infty}|g(r)| = +\infty, \qquad \lim_{r\to+\infty}\frac{f'(r)}{g'(r)} =: L \in \{-\infty\} \cup \mathbb{R} \cup \{+\infty\}, \tag{8}$$

*то существует предел*

$$\lim_{r\to+\infty}\frac{f(r)}{g(r)} = L = \lim_{r\to+\infty}\frac{f'(r)}{g'(r)}. \tag{9}$$

Специфика приведённого варианта правила Лопиталя в том, что в отличие от его традиционных форм здесь в условиях (8) не требуется традиционного условия существования бесконечного предела $\lim_{r\to+\infty}|f(r)| = +\infty$ для функции $f$ из числителей в (8) и (9).

Применим правило Лопиталя к функциям $g \coloneqq \ln M$ и $f \coloneqq \ln V$ с $L \coloneqq \rho$ из (7). Тогда по свойствам модельной функции роста $M$ из Определения в обозначении (2) существует

$$\lim_{r\to+\infty} \boldsymbol{\rho}_M(r) = \lim_{r\to+\infty}\frac{\ln V(r)}{\ln M(r)} = \lim_{r\to+\infty}\frac{f(r)}{g(r)} = L = \rho.$$

Таким образом, доказано существование предела (3) из утверждения II, а также установлено равенство (5). Наконец, из равенств (6) и (5) получаем

$$\lim_{r\to+\infty}\frac{M(r)}{M'(r)}\boldsymbol{\rho}'_M(r)\ln M(r) = \lim_{r\to+\infty}\boldsymbol{\rho}_M(r) - \lim_{r\to+\infty}\frac{M(r)V'(r)}{M'(r)V(r)} = \rho - \rho = 0,$$

что даёт существование предела (4), равного нулю, и завершает доказательство импликации I ⇒ II. Теорема доказана.

Напомним классическое определение уточнённого порядка. Дифференцируемая функция $\boldsymbol{\rho} \geq 0$, определённая на луче $R_{\to} \subset \mathbb{R}^+$, называется *уточнённым порядком в смысле Валирона* [1]–[5], если существуют два конечных предела

$$\rho \coloneqq \lim_{r\to+\infty}\boldsymbol{\rho}(r) \in \mathbb{R}^+, \qquad \lim_{r\to+\infty} r\boldsymbol{\rho}'(r)\ln r = 0. \tag{10}$$

**Следствие.** *Пусть функция $\boldsymbol{\rho} \geq 0$ дифференцируема на луче $R_{\to} \subset \mathbb{R}^+$. Тогда эквивалентны следующие два утверждения:*

I. *$\boldsymbol{\rho}$ – уточнённый порядок в смысле Валирона и $\rho \coloneqq \lim_{r\to+\infty}\boldsymbol{\rho}(r) \in \mathbb{R}^+$, как в первом равенстве из* (10).

II. *Функция $V(r) \coloneqq r^{\boldsymbol{\rho}(r)}$, $r \in R_{\to}$, – уточнённая функция роста относительно тождественной модельной функции роста $M \coloneqq \mathrm{id}\colon r \mapsto r$, что означает существование хотя бы одного предела из* (1), *а значит и каждого предела в* (1):

$$\rho = \lim_{r\to+\infty}\frac{rV'(r)}{V(r)} = \lim_{r\to+\infty}\frac{(\ln V(r))'}{(\ln r)'} = \lim_{x\to+\infty}\frac{(\ln v(x))'}{(\ln e^x)'} = \lim_{x\to+\infty}\frac{v'(x)}{v(x)} \in \mathbb{R}^+, \tag{11}$$

*где, как и прежде в Определении, $v(x) := V(e^x)$, $x \in \ln R_{\to} := \{\ln r : r \in R_{\to}\}$.*

*Доказательство.* Утверждение I Следствия – это, ввиду определения уточнённого порядка (10), в точности утверждение II Теоремы при $M \coloneqq \mathrm{id}$, где соотношения (3) и (4) – это соответственно первое и второе соотношения в (10). Утверждение II Следствия совпадает в случае $M \coloneqq \mathrm{id}$ с утверждением I Теоремы. Величина $\rho$ в (11) совпадает с величиной $\rho := \lim\limits_{r\to+\infty} \boldsymbol{\rho}(r) \in \mathbb{R}^+$ из утверждения I Следствия по равенству (5) Теоремы. Тем самым, Следствие доказано.

**Примеры.** Опишем широкий класс модельных функций роста, тесно связанных с важнейшими характеристиками роста субгармонических функций. Пусть $u$ – произвольная субгармоническая функция на $\mathbb{C}$, отличная от тождественной постоянной; $r \in \mathbb{R}_*^+$,

$$\boldsymbol{C}_u(r) := \frac{1}{2\pi}\int_0^{2\pi} u\left(re^{it}\right) dt, \qquad \boldsymbol{B}_u(r) := \frac{2}{\pi r^2}\int_0^r \boldsymbol{C}_u(s) s ds\ , \quad \boldsymbol{M}_u(r) := \sup\{u(z): |z| = r\}$$

– соответственно *интегральное среднее по окружности и по кругу с центром в нуле радиуса* $r$*, а также точная верхняя грань на этой окружности функции* $u$ [9; определение 2.6.7]. Тогда эти средние $\boldsymbol{C}_u$, $\boldsymbol{B}_u$ и верхняя грань $\boldsymbol{M}_u$ – модельные функции роста при условии их дифференцируемости при больших $r \in \mathbb{R}_*^+$ [9; теорема 2.6.8].

**Основная теорема.** *Пусть* $M$ *– модельная функция роста, а возрастающая строго положительная функция* $A$ *конечного порядка относительно* $M$ *в том смысле, что*

$$\limsup_{r\to+\infty} \frac{\ln\left(1 + A(r)\right)}{\ln M\,(r)} < +\infty.$$

*Тогда существует уточнённая функция роста* $V$ *относительно модельной функции роста* $M$*, с которой выполнено предельное соотношение*

$$\limsup_{r\to+\infty} \frac{A(r)}{V(r)} = 1.$$

В случае классического уточнённого порядка этот результат давно известен [1], [2; гл. I, § 12, теорема 16], [3; гл. II, § 2, теорема 2.1], [10; § 2, теорема 6] как основной и мотивировавший введение и использование понятия уточнённого порядка в теории роста классов функций. Уже в этом классическом случае доказательство технически довольно трудоёмко. Мы предполагаем привести полное доказательство Основной теоремы позже в ином месте.

## Литература

1. Valiron G. Lecture on the General Theory of Integral Functions. Toulouse, 1923. 234 p.
2. Левин Б. Я. Распределение корней целых функций. М.: ГИТТЛ, 1956. С. 632.
3. Гольдберг А. А., Островский И.В. Распределение значений мероморфных функций. М.: Наука, 1970. С. 591.
4. Bingham N. H., Goldie C. M., Teugels J. L. Regular variation. Encyclopedia Math. Appl., 27. Cambrige: Cambrige University Press, 1987. 494 p.
5. Гришин А. Ф., Малютина Т. И. Об уточнённом порядке // Комплексный анализ и математическая физика. Красноярск: Красноярский госуниверситет. 1998. С. 10–24.
6. Hörmander L. Notions of Convexity. Progress in Mathematics. Boston: Birkhäuser, 1994. 416 p.

7. Kiselman Ch. O. Order and type as measures of growth for convex or entire functions // Proceedings of the London Mathematical Society (3). 1993. Vol. 66. Pp. 152–186.

8. Taylor A. E. L'Hospital's Rule // American Mathematical Monthly. 1952. Vol. 59, No. 1. Pp. 20–24.

9. Ransford Th. Potential Theory in the Complex Plane. Cambridge: Cambrige University Press, 1995. 232 p.

10. Гришин А. Ф., Поединцева И. В. Абелевы и тауберовы теоремы для интегралов // Алгебра и анализ. 2014. Т. 26. № 3. С. 1–88; English Transl.: *St. Petersburg Math. J.*, **26**:3 (2015), 357–409.

## Variations on the proximate order

B. N. Khabibullin

*Bashkir State University*
*32 Zaki Validi Street, 450076 Ufa, Republic of Bashkortostan, Russia.*

*Email: khabib-bulat@mail.ru*

The concept of proximate order is widely used in the theories of entire, meromorphic, subharmonic and plurisubharmonic functions. We give a general interpretation of this concept as a proximate growth function relative to a model growth function. If a function $V$ is the proximate growth function with respect to the identity function on the positive semi-axis, then the function $\ln V$ is the classical proximate order. Our definition uses only one condition. This form of definition is also new for the classical proximate order.